# Dynamic charging management for electric vehicle demand responsive transport

Tai-Yu Ma

Luxembourg Institute of Socio-Economic Research (LISER), 11 Porte des Sciences, 4366 Esch-sur-Alzette, Luxembourg

Email: tai-yu.ma@liser.lu
ORCID: 0000-0001-6900-098X

## Abstract

With the climate change challenges, transport network companies started to electrify their fleet to reduce CO2 emissions. However, such ecological transition brings new research challenges for dynamic electric fleet charging management under uncertainty. In this study, we address the dynamic charging scheduling management of shared ride-hailing services with public charging stations. A two-stage charging scheduling optimization approach under a rolling horizon framework is proposed to minimize the overall charging operational costs of the fleet, including vehicles' access times, charging times, and waiting times, by anticipating future public charging station availability. The charging station occupancy prediction is based on a hybrid LSTM (Long short-term memory) network approach and integrated into the proposed online vehicle-charger assignment. The proposed methodology is applied to a realistic simulation study in the city of Dundee, UK. The numerical studies show that the proposed approach can reduce the total charging waiting times of the fleet by 48.3% and the total charged amount of energy of the fleet by 35.3% compared to a need-based charging reference policy.

## 1. Introduction

User-centered shared mobility services have been widely studied to reduce personal car use and enhance the accessibility of transit services in low-demand areas [3, 12]. With the climate change crisis, many transport network companies (TNC) started deploying electric vehicles for ride-hailing services to reduce CO2 emissions. However, due to limited driving range and long charging time, how to minimize electric vehicle (EV) charging costs and times given stochastic customer demand and uncertain public charging station availability has become an active research issue in recent years [5, 9, 13].

While the charging scheduling problems have been widely studied in the con-text of logistics, these studies focus mainly on deterministic green vehicle routing problems: partial recharge, non-linear charging function, capacitated charging station consideration [2, 6]. However, electrified ride-hailing service operational policy design presents a more challenging issue that needs to jointly optimize the fleet management and charging scheduling under uncertainty. The problem involves the decisions related to determining when and how much energy to charge by considering vehicles' driving needs and charged energy costs, and where to assign vehicles to charge during the daytime to minimize the waiting time and charging time with uncertain availability of public fast/rapid charging stations. A recent empirical study shows that electric ride-hailing service operations mainly rely on public DC fast-charging stations to save charging time [4]. As the number of DC fast chargers in a city is limited due to its high investment cost, an efficient charging management strategy needs to be developed by considering stochastic waiting time [7]. However, existing studies either assume uncapacitated charging stations or do not consider the stochastic charging demand of other EVs. Furthermore, most studies assume constant energy prices for recharging, while considering varying energy prices in vehicle routing-related problems are still limited [8]. While the location planning of the charging station might affect the operational costs of the operator, it is out of the scope of this study.

In this work, a dynamic and predictive charging scheduling and vehicle-charging station assignment strategy with public charging stations are proposed for dynamic shared ride-hailing services using a fleet of EVs. The novelty of this research is considering the varying energy prices and uncertain public charging station availability to minimize the overall charging costs and waiting times at charging stations of the fleet while satisfying customers' demands. The main contributions are summarized as follows.

a. Incorporate time-of-use energy price for vehicle charging scheduling optimization to minimize the overall charging costs of the fleet.
b. Integrate a predictive model using hybrid long short-term memory (LSTM) neural networks for public charging station occupancy prediction [10] into online vehicle-charging station assignment to minimize the overall charging operational time of the fleet. The results show that this predictive assignment strategy could reduce vehicle waiting times and energy costs significantly.
c. A numerical study is conducted using public charging session data of the city of Dundee, UK, to evaluate the performance of the proposed methodology.

## 2. Methodology

Consider a TNC operating a fleet of EVs for providing shared ride-hailing service. Customers arrive stochastically on short notice with requested pick-up and drop-off locations with desired pickup times. A vehicle dispatching and routing policy (e.g. [1, 11]) is applied to provide customers with door-to-door mobility service. The fleet of vehicles is charged to full at the beginning of the day and then recharge during the daytime based on an optimized day-ahead charging plan (described later). We assume that the maximum allowed recharging level is 80% or 90% of battery capacity and the charging rate is linear. Different from existing studies that assume uncapacitated operator-owned charging infrastructure, we consider the problem of charging scheduling using public rapid charging stations with stochastic charging demand from other EVs. In doing so, the waiting time when a vehicle arrives at a charging station is stochastic, depending on the day-to-day public charging station occupancy in the service area.

We propose a vehicle battery replenishment model to minimize the overall charging costs of the fleet, i.e. a set of EVs, $K = \{1, ..., \bar{k}\}$. Let $H = \{1, ..., h, ..., \bar{h}\}$ denote the charging decision planning horizon, discretized into $\bar{h}$ charging decision epochs with a homogenous time interval $\Delta$ (e.g. 30 minutes). The energy price on epoch $h$ is denoted as $\vartheta_h$, unchanged within the same epoch but might vary over the different epochs. We formulate this charging scheduling problem as a mixed-integer linear programming (MILP) below. The objective function (1) minimizes the total charging costs of the fleet over the planning horizon, where $y_{kh}$ is a binary decision variable being 1 if the vehicle k is charged on epoch $h$. $u_{kh}$ denotes the decision of the amount of energy to be charged on epoch $h$ for vehicle k. $\vartheta_h$ is the energy price on epoch $h$. $\bar{c}$ is a fixed cost for recharge, estimated as the average energy consumption costs to reach the charging stations per recharge operation. Constraint (2) is the energy conservation of vehicles with $e_{kh}$ being the state of charge (SOC) of vehicle k at the beginning of epoch h. $d_{kh}$ is the average energy consumption of vehicle $k$ on epoch $h$, obtained from vehicles' historical driving patterns. Constraint (3) assures that the SOC of vehicle k after recharging on epoch $h$ is no less than a reserve energy level $e_{min}$ plus the expected energy consumption on that epoch.

Constraint (4) ensures that $u_{kh}$ can be positive when $y_{kh} = 1$. Constraints (5)-(7) setup the range of vehicles' SOC and the maximum amount of energy $u_{max}$ that can be charged on one epoch, depending on the type of chargers used for vehicles.

$$(P1) \quad \min \sum_{k=1}^{\bar{k}} \sum_{h=1}^{\bar{h}} [\vartheta_h u_{kh} + \bar{c} y_{kh}] \quad (1)$$

subject to

$$e_{k,h+1} = e_{kh} + u_{kh} - d_{kh}, \quad \forall k \in K, h \in H \quad (2)$$

$$e_{kh} + u_{kh} \geq d_{kh} + e_{min}, \quad \forall k \in K, h \in H \quad (3)$$

$$u_{kh} \leq M y_{kh}, \quad \forall k \in K, h \in H \quad (4)$$

$$e_{k0} = e_{init}^k \quad (5)$$

$$e_{min} \leq e_{kh} \leq e_{max}, \quad \forall k \in K, h \in H \quad (6)$$

$$0 \leq u_{kh} \leq u_{max}, \quad \forall k \in K, h \in H \quad (7)$$

$$y_{kh} \in \{0,1\} \, \forall k \in K, h \in H \quad (8)$$

The above charging scheduling model provides an approximate charging plan for the fleet based on historical driving patterns of vehicles, disregarding different uncertain factors related to charging station capacity constraints and vehicle access distance to charging stations, and stochastic availability of public charging stations. These factors affect the vehicles' waiting time and charging operational costs, which will be optimized based on an online vehicle-charger assignment model below.

An online vehicle-charger assignment model extended from [9] is proposed by considering public charging stations for which other EVs might compete with the available charging resource (arrival time and charging time of other EVs are unknown). For each charging decision epoch $h$, we solve the following online vehicle-charger assignment problem to determine where to recharge for the scheduled-to-charge vehicles obtained from the day-ahead charging scheduling model of P1 (eqs. (1)-(8)). The predictive dynamic charging scheduling and vehicle-charger assignment framework is shown in Figure 1. The input data includes historical vehicle driving patterns of the fleet, experienced charging waiting times at the charging stations, locations of the charging stations, and energy prices in the study area. First, problem P1 is solved one day ahead while problem P2 (eqs. (9)-(18) below) is solved at the beginning of each decision epoch $h \in H$. Different from [9], we integrate the charging station occupancy prediction (on individual charger level) based on the hybrid LSTM neural networks (hybrid LSTM block in Figure 1) under the rolling horizon framework 1)(i.e. next hour from the beginning of decision epoch h) [10]. At the end of the day, the operator updates vehicle experienced waiting times and driving patterns as input for the P1 problem.

This P2 problem is formulated as a MILP as (9)-(18) by extending the model of [9].

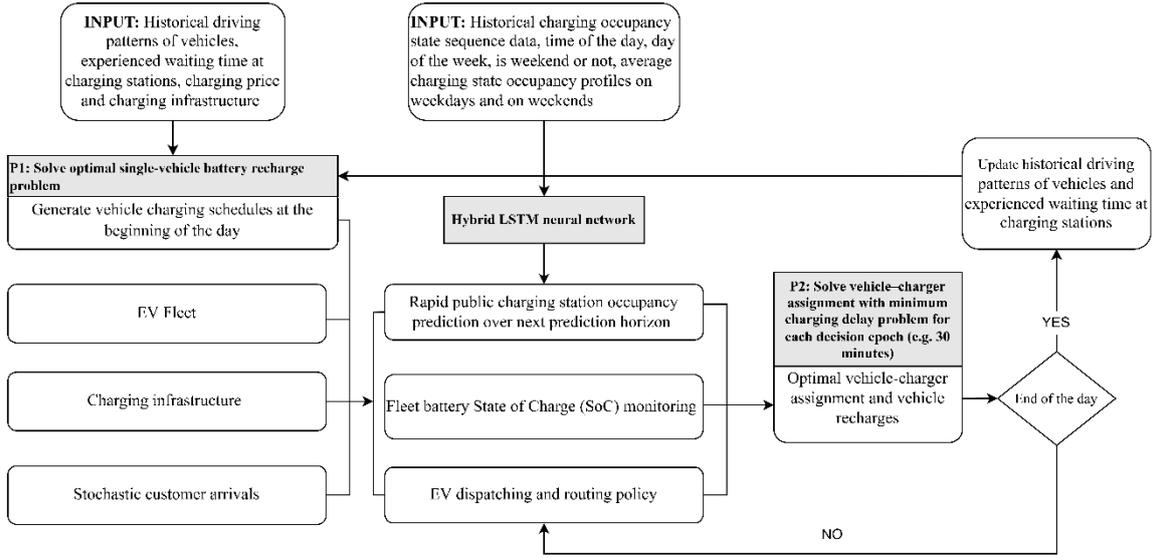

**Fig. 1.** Predictive dynamic charging scheduling and assignment framework for shared ride-hailing service.

*Notation*

| | |
|---|---|
| $I$ | Set of vehicles |
| $J$ | Set of chargers |
| $t_{ij}$ | Travel time from the location of vehicle $i$ to that of charger $j$ |
| $d_{ij}$ | Travel distance from the location of vehicle $i$ to that of charger $j$ |
| $e_i$ | Energy level of vehicle $i$ at the beginning of epoch $h$ (index $h$ is dropped) |
| $e_i^*$ | Energy level of vehicle $i$ after recharge at the end of epoch h, determined by the charging schedule obtained from solving the P1 problem (index $h$ is dropped) |
| $a_{ij}, b_{ij}$ | Predicted start and end time of a charging session of charger $j$ when vehicle $i$ arrives at charger $j$. Let $(g_{j1}, g_{j2}, \ldots, g_{js})$ be the predicted sequence of start/end charging session time series over a predictive horizon (e.g. one-hour ahead) with $g_{j1}$ and $g_{js}$ being the start time and end time of the prediction horizon, respectively. For a time slot $r$, $[g_{jr}, g_{j,r+1}]$ with $g_{jr} \leq t_{ij} \leq g_{j,r+1}$, if charger j is not available, then $a_{ij} = g_{jr}$, $b_{ij} = g_{j,r+1}$, otherwise $a_{ij} = b_{ij} = g_{js}$ |
| $\mu$ | Energy consumption rate of vehicles (kWh/km) |
| $\varphi$ | Charging rate of chargers (kW/min.) |
| $M_1, M_2$ | Large positive numbers |
| $\theta_1, \theta_2$ | Parameters |

Decision variable

| | |
|---|---|
| $X_{ij}$ | Vehicle $i$ is assigned to charger $j$ for a recharge if $X_{ij} = 1$, and 0 otherwise |
| $Y_{ij}$ | Amount of energy recharged at charger $j$ for vehicle $i$ |
| $S_{ij}$ | An artificial variable representing the waiting time of vehicle $i$ at charger $j$ |
| $W_j$ | Binary variable |

$$(P2) \min Z = \sum_{i \in I} \sum_{j \in J} t_{ij} X_{ij} + \theta_1 \sum_{i \in I} \sum_{j \in J} Y_{ij}/\varphi_j + \theta_2 \sum_{i \in I} \sum_{j \in J} S_{ij} \quad (9)$$

subject to

$$\sum_{j \in J} X_{ij} = 1, \quad \forall i \in I \quad (10)$$

$$\sum_{i \in I} X_{ij} \leq 1, \quad \forall j \in J \quad (11)$$

$$e_{min} \leq e_i - \mu d_{ij} X_{ij} + M_1(1 - X_{ij}), \forall i \in I, j \in J \quad (12)$$

$$e_i^* \leq Y_{ij} + e_i - \mu d_{ij} X_{ij} + M_1(1 - X_{ij}), \forall i \in I, j \in J \quad (13)$$

$$Y_{ij} \leq M_1 X_{ij}, \forall i \in I, j \in J \quad (14)$$

$$a_{ij} - t_{ij} - M_1(1 - X_{ij}) \leq M_2(1 - W_j), \forall i \in I, j \in J \quad (15)$$

$$b_{ij} - t_{ij} X_{ij} - M_1(1 - X_{ij}) \leq S_{ij} + M_2(1 - W_j), \forall i \in I, j \in J \quad (16)$$

$$a_{ij} - t_{ij} + M_1(1 - X_{ij}) + M_2 W_j > 0, \forall i \in I, j \in J \quad (17)$$

$$X_{ij}, W_j \in \{0,1\}, Y_{ij}, S_{ij} \geq 0, \forall i \in I, j \in J \quad (18)$$

The objective function (9) minimizes the total charging operational time as a weighted sum of the vehicle's access time to charging stations $t_{ij}$, charging times $Y_{ij}/\varphi_j$ and expected waiting time at charging stations $S_{ij}$. Constraints (10) and (11) state that each vehicle is assigned to exactly one charger and each charger can be assigned at most one vehicle. Constraint (12) states the SOC of vehicles is no less than the reserve level $e_{min}$ when arriving at a charger. Constraint (13) states that the SOC of vehicles after recharge needs to be no less than the planned battery level $e_i^*$. Constraint (14) ensures that the amount of energy to be recharged is positive when there is a charging event. Constraint (15)-(17) calculates the vehicle's waiting times when arriving at charging stations. Note that in case that $|I| \geq |J|$, constraints (10) and (11) need to be revised accordingly [9].

Figure 2 explains how the expected waiting time is calculated by eqs. (15)-(17). In Figure 2(b) if a charger $j$ is predicted to be occupied by another vehicle within $[g_{j3}, g_{j4}]$, vehicle $i$'s expected waiting time when arriving at $t_{ij}$ will be $g_{j3} - t_{ij}$. The charging occupancy prediction is based on a hybrid LSTM model by incorporating relevant features including time of day, day of the week, whether the day is weekday/weekend, average charging occupancy rate profiles on weekdays and weekends, and historical k-step backward charging occupancy states. The reader is referred to [10] for a more detailed description. The P1 problem can be easily solved using commercial solvers. The P2 problem is a generalized assignment problem for which the Lagrangian relaxation algorithm developed in [9] is applied to obtain near-optimum solutions for large test instances.

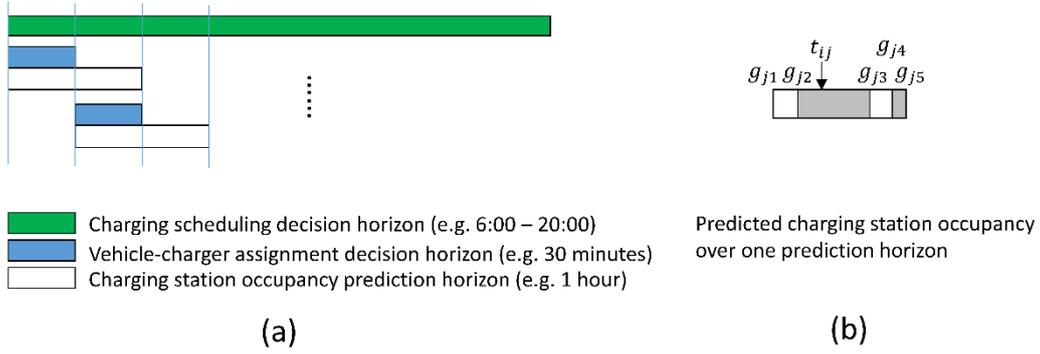

**Fig. 2.** (a) Concept framework of different decision/prediction horizons; (b) Waiting time estimation based on the predicted charging station occupancy using hybrid LSTM approach [10].

## 3. Numerical study

### 3.1 Study area and charging session data

We consider a simulation case study of shared ride-hailing service in the city of Dundee, UK (Figure 3). Customers' ride requests are generated randomly in the study area. A TNC operates a fleet of small electric capacitated vehicles to provide users door-to-door microtransit service. Vehicles' dispatching and routing policy is based on a non-myopic policy by considering future system costs when inserting a new customer on existing routes of vehicles [11]. Customer demand for a full day is assumed 800, randomly generated in the study area following some empirical probability distribution of their arrival time. The fleet size is assumed as 40 4-seater electric vehicles. A depot is located around the center of the study area. We assume that the fleet is fully charged at the beginning of the day and due to the battery capacity constraint, vehicles need to be recharged in the daytime using public rapid chargers to ensure their SOCs are always no less than a reserve level (20% of the battery capacity). We assume that the TNC charges their vehicles using rapid chargers only [4]. The public charging session data (https://data.dundeecity.gov.uk/dataset/ev-charging-data) in the study area is used for the charging station occupancy prediction [10]. This prediction model for a one-hour ahead prediction has an average accuracy of 81.87%. An illustrative example of the observed and predicted occupancies on a rapid charger is shown in Figure 4. We can observe that the model predicts approximately the charging occupancy pattern of the rapid charger.

We generate randomly 15 demand datasets with 800 random ride requests corresponding to 15 weekdays (3 weeks from 14/5/2018 to 1/6/2018). Different from [9], the energy consumption needs $d_{kh}$ in equation (2) is obtained using the first two-week datasets with a need-based policy, i.e. a vehicle goes to recharge at a nearest rapid charger whenever its SOC is lower than 25% of the battery capacity. An example of the customer's arrival time distribution is shown in Figure 5. Table 1 reports the parameters used for the simulation studies.

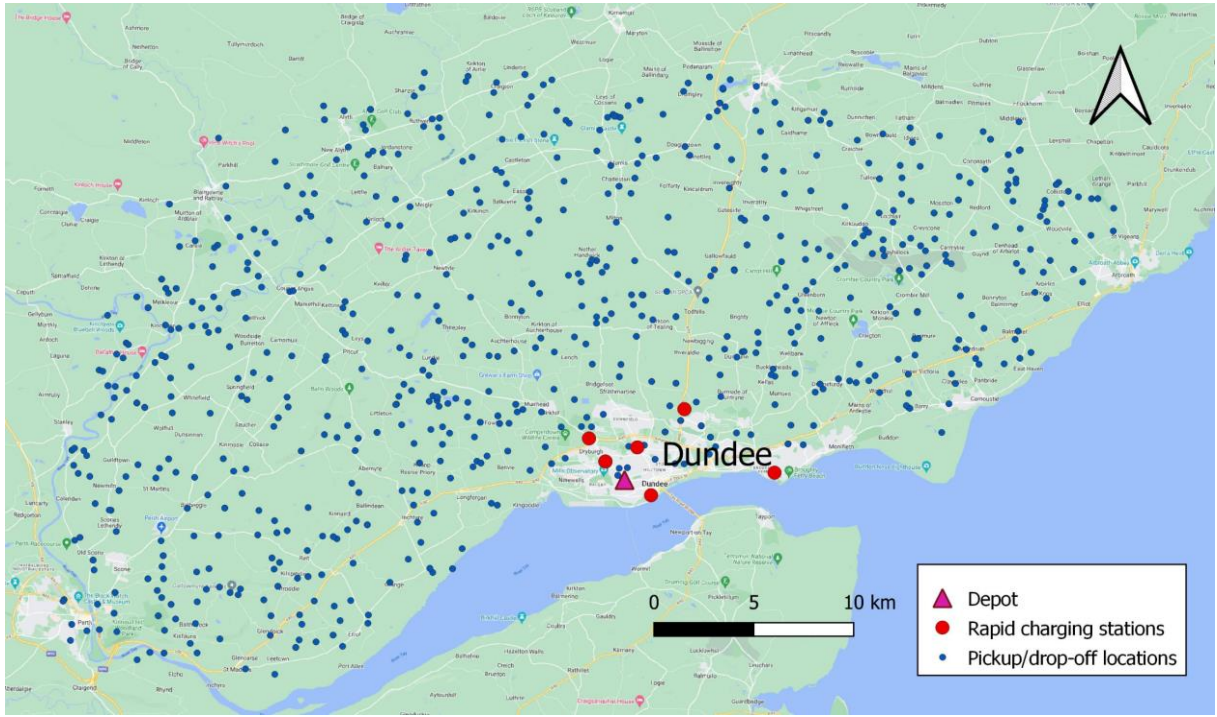

**Fig. 3.** Simulation case study on the city of Dundee, UK. The blue points are randomly generated pickup/drop-off locations of customers. There is one depot and 9 rapid chargers located in 6 different locations.

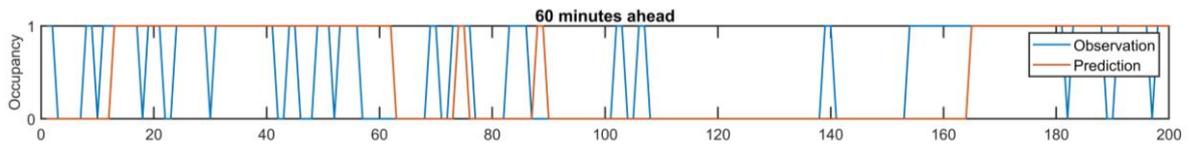

**Fig. 4.** Example of observed and predicted occupancies of a rapid charger on a weekday [10]. The y-axis is the binary state (occupied or not by an EV) of a charger. The x-axis is time (in minutes), starting from a reference point. The predicted occupancy is reported for 60 minutes ahead over time. We can observe that the prediction captures well the charging station occupancy patterns over different vehicle arrival intensities.

**Table 1.** The parameter setting for the simulation case study.

| | | | |
|---|---|---|---|
| Number of customers | 800 | $\vartheta$ | 0.25[2] (£/kWh) |
| Number of vehicle depots | 1 | $\mu$ | 0.204 (kWh/km) |
| Fleet size | 40 | $\bar{c}$ | 0.3 £ |
| Capacity of vehicles | 4 pers./veh. | $\Delta$ | 30 min. |
| Vehicle speed | 65 km/hour | $T$ | 6:30–22:00 |
| Battery capacity (B) | 62 kWh[1] | $\varphi_{DC\ fast}$ | 50/60 (kW/min.) |
| Number of rapid chargers | 9 | $\beta$ | 0.025 |
| $e_{min}$ | $0.2B$ | $\gamma$ | 0.5 |
| $e_{max}$ | $0.8B$ | | |

[1]Nissan LEAF E+. https://ev-database.org/car/1144/Nissan-Leaf-eplus
[2]https://www.dundeecity.gov.uk/news/article?article_ref=4137#:~:text=The%20new%20tariffs%20will%20be,when%20Scottish%20Government%20subsidies%20ended.

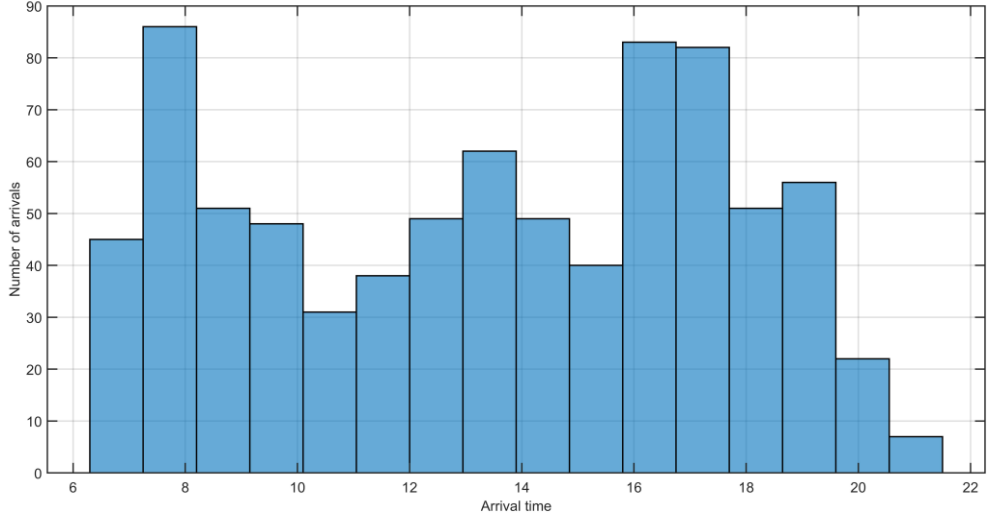

**Fig. 5.** Example of customer arrival times in the study area.

The discrete event simulation technique is applied to include queueing delays at charging stations. The simulation pseudocode is described in Algorithm 1.

**Algorithm 1. Simulation pseudocode for optimal vehicle-charging station assignment.**

| | |
|---|---|
| 1: | **Input**: Initial locations of vehicles for $I = \{1,2,\ldots,|I|\}$, charging plans $\{e^*_{ih}\}_{h\in H}$, $\forall i \in I$ (obtained from solving the P1 problem), locations and characteristics of chargers $J = \{1,2,\ldots,|J|\}$, and stochastic customer demand (pickup and drop-off locations of customers). |
| 2: | Upon the arrival of a new event E at time $t$ |
| 3: |   **Switch** Event_Type(E) |
| 4: |     // vehicle ($i$) arrives at its assigned charging station ($j$) |
| 5: |     **Case ARRIVAL_AT_CHARGER**: |
| 6: |       **If** no other EVs at charger $j$, **then** |
| 7: |         Charge immediately, waiting time $t^W_i$=0. |
| 8: |       **Else** |
| 9: |         Compute waiting time $t^W_i$ of the vehicle based on the duration that the other vehicles in the queue finish their charging. |
| 10: |       **End** |
| 11: |       Compute the vehicle's charging time as $t^G_i = [e^*_{ih(t)} - e_i(t)]/\varphi_j$. Schedule a new 'LEAVE_CHARGER' event at the time $t + t^W_i + t^G_i$. Update the vehicle's location and the state of the vehicle as 'CHARGING'. |
| 12: |     // vehicle leaves a charger after finishing a charging session |
| 13: |     **Case LEAVE_CHARGER**: update the state of the vehicle ($i$) to "AVAILABLE". Update the state of the charging station ($j$) and the energy state of the vehicle after recharge. |
| 14: |     **Case ARRIVAL_NEW_REQUEST**: Update vehicle states (locations, passengers on board, and energy level) until $t$. |
| 15: |     **If** a vehicle has no passengers on board and is marked as "go-charging", **Then** Schedule a new "ARRIVAL_AT_CHARGER" event at the time $t + t_{ij}$, where $t_{ij}$ is the travel time from the location of vehicle $i$ to the location of its assigned charger $j$. |
| 16: |     **Otherwise** |
| 17: |     Assign the new request to a vehicle based on the vehicle routing policy [11].     **End** |
| 18: |   **End Switch** |

| 19: | // charging station assignments |
|---|---|
| 20: | **If** $t \geq h\Delta$, **then** |
| 21: | Get the list of vehicles to be recharged in epoch $h$ according to all vehicles' charging plans $e^*$. If the number of vehicles to be recharged is greater than the number of chargers, postpone the additional vehicles with higher battery levels (lower priority) to the next charging epoch $h + 1$. The resulting list of vehicles to be recharged is denoted as list_recharge. |
| 22: | Solve the P2 problem of eqs. (9)-(18) for the vehicles on list_recharge and obtain the vehicle-charger assignment pair. |
| 23: | **For** $i \in$ list_recharge, If vehicle *i* is idled, schedule a new "ARRIVAL_AT_CHARGER" event at the time $t + t_{ij}$. Mark the status of the vehicle *i* as "go-charging". |
| 24: | **End** |
| 25: | $h \leftarrow h + 1$ |
| 26: | **End** |
| 27: | Continue until all scheduled events are executed. |
| 28: | Every vehicle continues its unfinished tour until all customers are served and then returns to its initial depot. |

## 1.1 Results

To validate the proposed methodology, we compare the performance of the proposed non-myopic optimal charging policy (OCP)(eqs. (1)-(18)) with the need-based policy. We evaluate also the benefits of with (OCP1) or without(OCP0) predictive information of future charging occupancy states in terms of charging waiting time savings. The simulation results are the average performance obtained from 5 test datasets (i.e. one represents a randomly generated 800 customers in the study area) during the third week.

Table 2 presents the results of the comparison. The performance measures in terms of average waiting time, average charging time, and average operational time for a charge are based on each charging session. The total waiting times for recharge of the fleet, total charging time of the fleet, and the total charged amount of energy of the fleet are presented in the last three columns. We can observe that using the OCP0 policy could reduce the average charging waiting time from 24.2 minutes to 20 minutes while applying OCP1 could lead to a significant waiting time reduction to 12.2 minutes (-49.5% compared to the need-based policy). When comparing the benefits of incorporating the predictive information, applying the OCP1 policy could reduce the average charging waiting time by -38.8%. In terms of the total charged amount of energy, the OCP1 policy could lead to significant charging costs savings by 35.3% compared to the need-based policy.

**Table 2.** Comparison of the need-based policy and the non-myopic optimal charging policy with and without predictive information.

| Charging policy[1] | ACWT[2] (min.) | ACT (min.) | AOTC (min.) | TWT (hour) | TCT (hour) | TCE (kWh) |
|---|---|---|---|---|---|---|
| NP | 24.2 | 43.4 | 69.4 | 35.7 | 63.3 | 3165.5 |
| OCP0 | 20.0 | 27.3 | 49.1 | 28.9 | 39.5 | 1976.9 |
| OCP0 vs. NS | -17.5% | -37.0% | -29.3% | -18.9% | -37.5% | -37.5% |
| OCP* | 12.2 | 27.2 | 41.6 | 18.4 | 41.0 | 2048.9 |
| **OCP* vs. NP** | **-49.5%** | **-37.2%** | **-40.1%** | **-48.3%** | **-35.3%** | **-35.3%** |
| OCP* vs. OCP0 | -38.8% | -0.4% | -15.2%% | -36.3% | 3.6% | 3.6% |

Remark: 1. NP: need-based policy; OCP0: Optimal charging policy **without** predictive information; OCP*: Optimal charging policy **with** predictive information. 2. **ACWT**: Average charging waiting time, **ACT**: Average charging time, **AOTC**: Average operational time for a recharge, **TWT**: Total waiting time of the fleet on a day, **TCT**: Total charging time of the fleet on a day, **TCE**: Total amount of charged energy of the fleet on a day.

## 2      Conclusions and discussion

Dynamic electric ridesharing fleet charging management under uncertainty is a challenging research issue. While existing studies mainly focus on deterministic electric vehicle routing problems in the logistic context, this study considers dynamic charging scheduling optimization problems with stochastic waiting time at public charging stations. A two-stage optimization approach is proposed to solve this problem and compared to a reference need-based charging policy. Different from our previous study [9], we integrate the predictive information of charging station occupancy in a rolling horizon framework to minimize the expected waiting time using public charging stations. The proposed methodology is tested using the public charging session data of the city of Dundee, UK. Our results show that the proposed charging strategy could lead to significant charging waiting time savings and reduce the amount of charged energy on the day, compared to the reference need-based charging policy.

The limitation of this research is that the day-to-day energy consumption of vehicles is stochastic and needs to be updated accordingly. Future extensions could consider a learning process to track vehicles' day-to-day energy consumption patterns to better estimate their energy needs and derive more adaptive charging plans under uncertain environments. Using a reinforcement learning approach or a hybrid approach to solve such dynamic decision-making problems under uncertainty are promising solutions to be studied in the future.

**Acknowledgement**

The work was supported by the Luxembourg National Research Fund (C20/SC/14703944).

**References**

1. Alonso-Mora, J. et al.: On-demand high-capacity ride-sharing via dynamic trip-vehicle assignment. Proc. Natl. Acad. Sci. 114, 3, 462–467 (2017). https://doi.org/10.1073/pnas.1611675114.
2. Asghari, M., Mirzapour Al-e-hashem, S.M.J.: Green vehicle routing problem: A state-of-the-art review. Int. J. Prod. Econ. 231, 107899 (2021). https://doi.org/10.1016/j.ijpe.2020.107899.
3. Chow, J. et al.: Spectrum of Public Transit Operations : From Fixed Route to Microtransit. (2020).
4. Jenn, A.: Electrifying Ride-sharing: Transitioning to a Cleaner Future, https://3rev.ucdavis.edu/policy-brief/electrifying-ride-sharing-transitioning-cleaner-future, last accessed 2021/12/08.
5. Keskin, M. et al.: Electric Vehicle Routing Problem with Time-Dependent Waiting Times at Recharging Stations. Comput. Oper. Res. 107, 77–94 (2019). https://doi.org/10.1016/j.cor.2019.02.014.
6. Keskin, M., Çatay, B.: Partial recharge strategies for the electric vehicle routing problem with time windows. Transp. Res. Part C Emerg. Technol. 65, 111–127 (2016). https://doi.org/10.1016/j.trc.2016.01.013.
7. Kim, J. et al.: Scheduling and performance analysis under a stochastic model for electric vehicle charging stations. Omega. 66, 278–289 (2017). https://doi.org/10.1016/j.omega.2015.11.010.
8. Klein, P.S., Schiffer, M.: Electric vehicle charge scheduling with flexible service operations. (2022). https://doi.org/https://arxiv.org/abs/2201.03972.

9. Ma, T.-Y.: Two-stage battery recharge scheduling and vehicle-charger assignment policy for dynamic electric dial-a-ride services. PLoS One. 16, 5, e0251582 (2021). https://doi.org/10.1371/journal.pone.0251582.
10. Ma, T.-Y., Faye, S.: Multistep electric vehicle charging station occupancy prediction using hybrid LSTM neural networks. Energy. 244, 123217 (2022). https://doi.org/10.1016/j.energy.2022.123217.
11. Ma, T.Y. et al.: A dynamic ridesharing dispatch and idle vehicle repositioning strategy with integrated transit transfers. Transp. Res. Part E Logist. Transp. Rev. 128, 417–442 (2019). https://doi.org/10.1016/j.tre.2019.07.002.
12. Ma, T.Y. et al.: A user-operator assignment game with heterogeneous user groups for empirical evaluation of a microtransit service in Luxembourg. Transp. A Transp. Sci. 17, 4, 946–973 (2021). https://doi.org/10.1080/23249935.2020.1820625.
13. Shi, J. et al.: Operating Electric Vehicle Fleet for Ride-Hailing Services with Reinforcement Learning. IEEE Trans. Intell. Transp. Syst. 21, 11, 4822–4834 (2020). https://doi.org/10.1109/TITS.2019.2947408.